\documentclass[a4paper,12pt]{book}
\usepackage{amsmath}
\usepackage[latin1]{inputenc}
\usepackage{amssymb}
\usepackage{graphicx}
\usepackage{makeidx}
\usepackage{mathrsfs}
\usepackage{amsfonts}
\usepackage[mathscr]{eucal}
\usepackage{amsmath}
\usepackage{fancyhdr}
\font\fmtitle=cmbx12 scaled \magstep2
\font\text=cmr10 at 12 truept
 at 11 truept
 at 10 truept

\begin{document}
 \thispagestyle{empty}
 
 \markboth
 {\it \centerline {Proceedings}}
 {\it \centerline{Rudolf GORENFLO}}  
\def\pni{\par \noindent}
\def\vsh{\vskip 0.25truecm\noindent}
\def\vs{\vskip 0.5truecm}
\def\vvs{\vskip 1.0truecm}
\def\vvvs{\vskip 1.5truecm}
\def\vsp{\vsh\par}
\def\vsn{\vsh\pni}
\def\cen{\centerline}
\def\ra{\item{(a)\ }} \def\rb{\item{(b)\ }}   \def\rc{\item{(c)\ }}
\def\eg{{e.g.}\ } \def\ie{{i.e.}\ }
\def\sg{\hbox{sign}\,}
\def\sgn{\hbox{sign}\,}
\def\sign{\hbox{sign}\,}
\def\e{\hbox{e}}
\def\exp{\hbox{exp}}
\def\ds{\displaystyle}
\def\dis{\displaystyle}
\def\q{\quad}    \def\qq{\qquad}
\def\lan{\langle}\def\ran{\rangle}
\def\l{\left} \def\r{\right}
\def\lra{\Longleftrightarrow}
\def\arg{{\rm arg}}
\def\argz{{\rm arg}\, z}
\def\argG{{x^2/ (4\,a\, t)}}
\def\d{\partial}
 \def\dr{\partial r}  \def\dt{\partial t}
\def\dx{\partial x}   \def\dy{\partial y}  \def\dz{\partial z}
\def\rec#1{\frac{1}{#1}}
\def\log{{\rm log}\,}
\def\erf{{\rm erf}\,}     \def\erfc{{\rm erfc}\,}

\def\NN{{\rm I\hskip-2pt N}}
\def\MM{{\rm I\hskip-2pt M}}
\def\RR{\vbox {\hbox to 8.9pt {I\hskip-2.1pt R\hfil}}\;}
\def\CC{{\rm C\hskip-4.8pt \vrule height 6pt width 12000sp\hskip 5pt}}
\def\II{{\rm I\hskip-2pt I}}
\def\erf{{\rm erf}\,}   \def\erfc{{\rm erfc}\,}
\def\exp{{\rm exp}\,} \def\e{{\rm e}}
\def\ss{{s}^{1/2}}   
\def\N{\bar N}  
\def\ss{{s}^{1/2}} 
\def\stt{{\sqrt t}}
\def\lst{{\lambda \,\stt}}
\def\Et{{E_{1/2}(\lst)}}
\def\u{\widetilde{u}}
\def\ul{\widetilde{u}} 
\def\uf{\widehat{u}} 
\def\bar{\widetilde}
\def\A{{\mathcal {A}}}
\def\L{{\cal L}} 
\def\F{{\cal F}} 
\def\M{{\cal M}}  
\def\Fdiv{\,\stackrel{{\cal F}} {\leftrightarrow}\,}
  \def\Ldiv{\,\stackrel{{\cal L}} {\leftrightarrow}\,}
  \def\Mdiv{\,\stackrel{{\cal M}} {\leftrightarrow}\,}
\def\args{(x/ \sqrt{a})\, s^{1/2}}
\def\argsa{(x/ \sqrt{a})\, s^{\beta}}
\def\barr{\widetilde}
\def\argG{ x^2/ (4\,a\, t)}
\def\G{{\cal {G}}}
\def\Gc{{\cal {G}}_c}	\def\Gcs{\barr{\Gc}} 
\def\Gs{{\cal {G}}_s}	\def\Gss{\barr{\Gs}} 
\def\f{\bar{f}}
\def\g{\bar{g}}
\def\u{\bar{u}}
\centerline{{\fmtitle Mittag-Leffler Waiting Time, Power Laws,}}
\vskip 0.20truecm  
\centerline{{\fmtitle Rarefaction, Continuous Time Random Walk,}}
\vskip 0.20truecm
\centerline{{\fmtitle Diffusion Limit}} %
\vskip 0.20truecm
\begin{center}
{{\bf Rudolf GORENFLO}}
\\
Free University Berlin, Germany
\\ E-mail: {\tt gorenflo@mi.fu-berlin.de}    
\end{center}
\centerline{{\bf Abstract}}
\noindent
We discuss some applications of the Mittag-Leffler function and related probability distributions in 
the theory of renewal processes and continuous time random walks. 
In particular we show the asymptotic (long time) equivalence of a generic power law waiting time to 
the Mittag-Leffler waiting time distribution via rescaling and respeeding the clock of time. 
By a second respeeding (by rescaling the spatial variable) we obtain the diffusion limit 
of the continuous time random walk under power law regimes in time and in space. 
Finally, we exhibit the time-fractional drift process as a diffusion limit of the fractional 
Poisson process and as a subordinator for space-time fractional diffusion.
\vskip .50truecm \noindent 
{AMS Subject Classification Numbers:}  26A33. 33E12, 45K05, 60G18,  60G50, 60G52, 60K05, 76R50.
\vskip .50truecm \noindent
{Keywords:} Continuous Time Random Walk, Fractional Diffusion, 
Mittag-Leffler Function, Power Laws, Rarefaction (Thinning), Renewal  Process, Subordination
\centerline{{\bf Contents}} 
1. Introduction
\\ 2. Continuous Time Random Walk (CTRW)
 \\ 3. Power Laws and Well-Scaled Passage to the Diffusion Limit
  \\ 4. Thinning (Rarefaction) of a Renewal Process under Power Law Regime
\\ 5. Mittag-Leffler Waiting Time and Space-Time-Fractional Diffusion
  \\ 6. Time-Fractional Drift and Subordination
  \\ 7. Conclusions
  \\  \phantom{8.} Acknowledgement
\\ \phantom{9.} References
\section*{1. Introduction}
 At begin of the past century Gösta Magnus Mittag-Leffler  introduced the entire function
$$  E_\alpha (z) := \sum_{n=0}^{\infty} \frac{z^n}{\Gamma (1 + \alpha n)}\,,
\quad z\in \CC\,,  \quad \Re \alpha >0\,.
\eqno (1.1)$$                                            
and investigated its basic properties, see Mittag-Leffler (1903). 
Although this function, named after him, and its generalization
$$ E_{\alpha,\beta} (z) := \sum_{n=0}^{\infty} \frac{z^n}{\Gamma (\beta + \alpha n)}\,,
\quad z\in \CC\,,  \q \Re \alpha >0\,,\;\Re \beta >0 \,.               \eqno (1.2) $$       
then was investigated by some authors, e.g. Wiman(1905), and used for the solution of the second kind Abel 
integral equation, see Hille and Tamarkin (1930), it did not find the deserved attention of the general community. 
In several important books and collections of formulas like Abramowitz and Stegun (1965) 
and Gradshteyn and Ryzhik (2000) on special functions it was ignored; 
a noteworthy exception is Chapter XVIII in Vol. III of Erdelyi et al.  (Bateman project of 1955). 
  
 In this report a prominent role as waiting time distribution will be played by the 
 {\it Mittag-Leffler probability distribution function}, see Pillai (1990),
$$   \Phi_\beta^{ML}(t)= 1 - E_\beta(-t^\beta)\,, \q t\ge 0\,,\q 0<\beta \le 1\,. \eqno(1.3)$$            
the probability of a waiting time greater than $t$, called {\it survival probability},                      
$$ \Psi_\beta^{ML} = 1-\Phi^{ML}_\beta(t) = E_\beta(-t^\beta)\,,$$
and the corresponding probability density function
$$\phi_\beta^{ML}(t) = \frac{d}{dt}\Phi_\beta^{ML}(t)= -\frac{d}{dt} E_\beta(-t^\beta)
= t^{\beta-1} \, E_{\beta,\beta}(-t^\beta)
\,, \; t\ge 0\,.
\eqno(1.4)$$
The functions $E_\beta(-t^\beta)$ and $\phi_\beta^{ML}(t)$ 
 are completely monotone (representable as Laplace transforms of non-negative measures, 
 see Feller (1971), 
 in concreto Gorenflo and Mainardi (1997),
$$
\phi_\beta^{ML}(t) = \frac{1}{\pi}\int_0^\infty \!\!
\frac{r^\beta \, \sin (\beta\pi)}{r^{2\beta} + 2r^\beta \, \cos(\beta \pi) +1}\, \exp(-rt)\, dr\,.
\eqno(1.5)$$                                  
Particularly important is the {\it power law} asymptotics for $t \to \infty$:
$$
E_\beta (-t^\beta)\sim \frac{t^{-\beta}}{\Gamma(1-\beta)}\,, \q
\phi_\beta^{ML}(t) \sim \frac {\Gamma(\beta+1)\, \sin(\beta \pi)}{\pi}\, t^{-\beta-1}
 \eqno(1.6)$$                                           
in contrast to the exponential decay of  $E_1(t)= \exp(-t)$.  
These asymptotics are an essential reason for the importance of these functions in modelling 
anomalous diffusion processes.  
 
      For later use let us  here note the Laplace transforms 
	  (for $0<\beta \le 1$  and $Re\, s \ge 0$)
 $$ \widetilde \Psi_\beta^{ML} = \frac{s^{\beta-1}}{s^\beta +1}\,, \q
\widetilde \phi_\beta^{ML} = \frac{1}{s^\beta +1}\,.    $$
\pni 
 {\bf Notations:} 
 $\widehat f(\kappa):= {\ds \int_{-\infty}^{+\infty}\!\! f(x)\,\exp(i\kappa x)\,dx}$,  
 with $\kappa$  real, 
 for the Fourier transform,
 $\widetilde g(s):= {\ds \int_0^\infty}\!\! g(t)\, \exp (-st)\, dt$, 
 with $s$ in a suitable right half-plane, for the Laplace transform.
\pni
   {\bf Remark:} Formula (1.5) exhibits for $0<\beta<1$ the Mittag-Leffler waiting time density as 
   a mixture of infinitely many exponential waiting time densities   
   $r\,\exp(-rt)$ with $r$-dependent 
   weight function behaving like $r^{-\beta-1}\sin(\beta\pi)$  for large $r$, 
   like  $r^{\beta-1}\sin(\beta\pi)$  for small  $r$. 
   Again we have power law asymptotics. 
   For interpretation consider the fact that the density   $r\, \exp(-rt)$
   whose mean is  $1/r$  decays exponentially fast for large $r$ but not so fast for 
   small $r$. 
   Starting around 1965 these functions attained increasing attention among researchers, 
   first in the theory of elasticity and relaxation (Caputo and Mainardi (1971), Nonnenmacher (1991)), 
   and later in the theory of continuous time random walk (pioneered 
   by Montroll and Weiss (1965) 
   who, however did not see the relevance of the Mittag-Leffler function) 
   and its limiting relation to fractional diffusion. 
   There are instances where researchers found the Laplace transform of solutions to 
   certain problems 
   but did not identify it as the transforms of functions of Mittag-Leffler type, 
   for example Gnedenko and Kovalenko (1968) in the theory of thinning or rarefaction 
   of a renewal process, Balakrishnan (1985) in his asymptotic investigation of 
   continuous time random walks. 
   For the latter it was Hilfer with Anton (1995) who clarified the relationship 
   between continuous time random walk (in the sense of Montroll and Weiss), 
   Mittag-Leffler waiting time and fractional derivative in time. 
   As more recent monographs with useful information on Mittag-Leffler functions 
   let us cite Samko, Kilbas and Marichev (1993), Podlubny (1999), Miller and Ross (1993), 
   Kilbas, Srivastava and Trujillo (2006), Mathai and Haubold (2008). 
   See also the comprehensive recent report by Haubold, Mathai and Saxena (2009). 
   Due to the growing importance of Mittag-Leffler functions there now is also activity 
   in the development of efficient methods for their numerical calculation, 
   see. e.g. Gorenflo, Loutschko and Luchko (2002) and Seybold and Hilfer (2008).
 \vsp   
	  In Section 2 of the present paper we will sketch the basic formalism 
	  of continuous time random walk, 
	  then in Section 3 under power law regime the well-scaled transition 
	  to the diffusion limit yielding the Cauchy problem for space-time fractional diffusion. 
	  Section 4 is devoted to thinning (or rarefaction) of a renewal process  
	  under power law regime and the relevant scaled transition via rescaling and respeeding, 
	  to a renewal process with Mittag-Leffler waiting time. 
	  Then, in Section 5, the Mittag-Leffler waiting time law and its 
	  relevance in continuous time random walk and the limiting fractional 
	  diffusion processes are discussed. 
	  Again, the transitions are achieved via re-scaling and re-speeding. 
	  Finally, in Section 6, we discuss the time-fractional drift and its role as a 
	  time-changing subordinator (producing the operational time from the physical time) 
	  in space-time fractional diffusion. Conclusions are drawn in Section 7.
  
\section*{2. Continuous Time Random Walk}
      Starting in the Sixties and Seventies of the past century the concept of continuous 
	  time random walk,  CTRW, became popular in physics as a rather general (microscopic) 
	  model for diffusion processes. 
	  Let us just cite Montroll and Weiss (1965), Montroll and Scher (1973), 
	  and the monograph of Weiss (1994). 
	  Mathematically, a CTRW is a {\it compound renewal process} or 
	  a {\it renewal process with rewards} or a 
	  {\it random-walk subordinated to a renewal process}, 
	  and has been treated as such by Cox (1967).  
	  It is generated by a sequence of independent identically distributed ($iid$)
	   positive waiting times $T_1,T_2,T_3, \dots$, 
	  each having the same probability distribution $\Phi(t)$, $t\ge0$, 
	  of waiting times $T$, and a sequence of $iid$ 
	  random jumps $X_1, X_2, X_3, \dots$, 
	  each having the same probability distribution function $W(x)$, $x \in  \RR$,
	  of jumps $X$. 
	  These distributions are assumed to be independent of each other. 
	  Allowing generalized functions (that are interpretable as measures) in the sense of 
	  Gelfand and Shilov (1964) we have corresponding probability densities
	  $\phi(t)=\Phi'(t)$  and $w(x)= W'(x)$ that we will use for ease of notation. 
	  Setting $t_0=0$, $t_n=T_1+T_2+\dots + T_n$  for $n\in \NN$, and
	  $x_0=0$, $x_n= X_1+X_2+\dots + X_n$, $x(t)=x_n$  
	  for $t_n\le t<t_{n+1}$  we get a (microscopic) model of a diffusion process. 
	  A wandering particle starts in  $x_0=0$ and makes a jump $X_n$ {\it at each instant} $t_n$.
	  Natural probabilistic reasoning then leads us to the {\it integral equation 
	  of continuous time random walk}  for the probability density $p(x,t)$  of 
	  the particle being in position $x$  at instant $t\ge 0$:
$$
 p(x,t) =  \delta (x)\, \Psi(t)\, +
  \int_0^t  \!\!  \phi(t-t') \, \left[
 \int_{-\infty}^{+\infty}\!\!  w(x-x')\, p(x',t')\, dx'\right]\,dt'\,,
\eqno(2.1)$$                
Here the {\it survival probability} 
$$\Psi(t) = \int_t^\infty \phi(t') \, dt'
\eqno(2.2)$$                                          
denotes the probability that at instant $t$ the particle still is sitting 
in its initial position $x_0=0$. 
Using, generically, for the Laplace transform of a function $f(t)$, $t\ge 0$, 
and the Fourier transform of a function $g(x)$, $x\in  \RR$, the notations
$$  \begin{cases} 
 \widetilde f(s)
 := {\ds \int_0^{\infty} \!\!\e^{\ds \, -st}\, f(t)\, dt \,,} \q \Re s \ge \Re s_0\,,\\ 
    \widehat g(\kappa)
  := {\ds \int_{-\infty}^{+\infty} \!\! \e^{\,\ds +i\kappa x}\,g(x)\, dx}
     \,, \; \kappa \in \RR\,,
	 \end{cases} 
\eqno(2.3)$$
we arrive, via   $\widehat \delta(\kappa) \equiv 1$
and the convolution  theorems, in the transform domain at the equation
$$
\widehat{\widetilde p}(\kappa ,s)
 =   \widetilde\Psi(s) +
   \widehat w(\kappa )\,\widetilde \phi(s) \,
   \widehat{\widetilde p}(\kappa ,s) \,, 
\eqno(2.4)$$        
which, by $\widetilde \Psi(s)= (1-\widetilde \phi(s))/s$
    implies the Montroll-Weiss equation, see Montroll and Weiss (1965),
	Weiss (1994),
$$\widehat{\widetilde p}(\kappa ,s) = 
\frac{1-\widetilde \phi(s))}{s} \, \frac{1}{1- \widehat w(\kappa )\,\widetilde \phi(s)}\,. 
\eqno(2.5)$$                       
Because of $|\widehat w(\kappa)| <1$,
$|\widetilde \phi(s)|<1$   for  $\kappa \ne 0$, $s \ne 0$,  we can expand into a geometric series
$$
\widehat{\widetilde p}(\kappa ,s) = \widetilde\Psi(s)\, \sum_{n=0}^{\infty}
   \,\left[\widetilde \phi(s) \,\widehat w(\kappa )\right]^n\,,
\eqno(2.6)$$                 
and promptly obtain the series representation of the CTRW, see Cox (1967), Weiss (1994),
$$
p(x,t)= \sum_{n=0}^{\infty} v_n(t)\, w_n(x)\,.
\eqno(2.7)$$                     
Here the functions  $v_n = (\Psi\,* \, \phi^{*n})$ and $w_n= w^*$ 
are obtained by iterated convolutions in time $t$ and in space $x$, respectively, 
in particular we have 
$$ v_0(t) =(\Psi \,*\, \delta)(t) = \Psi(t)\,, \;
v_1(t) =(\Psi \,*\, \phi)(t) \,, \; w_0(x)=\delta(x)\,,\; w_1(x)= w(x)\,.$$
The representation (2.7) can be found without the detour over (2.5) by direct probabilistic treatment. 
It exhibits the CTRW as a subordination of a random walk to a renewal process.    
  
   Note that in the special case $\phi(t)= m \,\exp(-mt)$, $m>0$,
    the equation (2.1) describes the compound Poisson process. 
	It reduces after some manipulations (best carried out in the transform domain) to
	 the {\it Kolmogorov-Feller equation}
     $$ \frac{\d}{\dt} p(x,t)= -m p(x,t) = m \int_{-\infty}^{+\infty} w(x-x')p(x',t)\, dx'$$                         
and from (2.7)  we obtain the series representation
      $$  p(x,t)=\e^{-mt}\, \sum_{n=0}^{\infty}\frac{(mt)^n}{n!} \, w_n(x)\,.$$  
	\section*{3. Power Laws and Well-Scaled Passage to the Diffusion Limit}
	  In recent decades  power laws in physical (and also economical and other) 
	  processes and situations have become increasingly popular for modelling slow 
	  (in contrast to fast, 
	  mostly exponential) decay at infinity. 
	  See Newman (2005) for a general introduction to this concept. 
	  For our purpose let us assume that the the distribution of jumps is symmetric, 
	  and that the distribution of jumps, likewise that of waiting times, 
	  either has finite second or first moment, respectively, or decays near infinity 
	  like a power with exponent $-\alpha$   or $-\beta$, respectively,
	  $ 0<\alpha<2$, $0<\beta<1$. 
	  Then we can state Lemma 1 and Lemma 2. 
	  These lemmata and more general ones (e.g. with slowly varying decorations 
	  of the power laws (a) and (b)) can be distilled from the Gnedenko theorem 
	  on the domains of attraction of stable probability laws (see Gnedenko and Kolmogorov (1954)), 
	  also from Chapter 9 of Bingham, Goldie and Teugels (1967). 
	  For wide generalizations (to several space dimensions and to anisostropy) see 
	  Meerschaert and Scheffler (2001). 
	  They can also be modified to cover the special case of 
	  smooth densities $w(x)$  and $\phi(t)$   
	  and to the case of fully discrete random walks, see Gorenflo and Abdel-Rehim (2004), 
	  Gorenflo and Vivoli (2003). 
	  For proofs see also Gorenflo and Mainardi (2009). 
\\
{\bf Lemma 1} (for the jump distribution):
\\ 
{\it Assume $W(x)$ increasing, $W(-\infty)=0\,,\, W(\infty)=1$, and symmetry  
$W(-x)+W(x)=1$
  for all continuity points $x$ of $W(x)$, and assume} (a) {\it or} (b):
\\
(a)     $\sigma^2:= {\ds \int_{-\infty}^{+\infty} x^2 \, dW(x)}<\infty$, labelled as
$\alpha=2$\,;
\\
 (b) ${\ds \int_x^\infty dW(x')}\sim b \alpha^{-1} x^{-\alpha}$ for $x\to \infty\,,
 \; 0<\alpha<2\,, \; b>0 \,.$
\\
{\it Then, with $\mu= \sigma^2/2$   in case} (a), 
{\it  $\mu = b\pi/[\Gamma(\alpha+1) \, \sin(\alpha \pi/2)]$
      in case} (b)
	  \\
{\it we have the asymptotics $1-\widehat w(\kappa)\sim \mu |\kappa|^\alpha$ for $\kappa \to 0$.}
\\
{\bf Lemma 2} (for the waiting time distribution): 
\\
{\it Assume  $\Phi(t)$ increasing, $\Phi(0)=0)$, $\Phi(\infty)=1$, and} (A) {\it or} (B). 
\\ (A)  {\it $\rho:={\ds \int_0^\infty\! t \, d\Phi(t)}<\infty$, labelled as $\beta=1$,}
\\ (B) {\it $1-\Phi(t)\sim c \beta^{-1}t^{-\beta}$   for $t\to \infty\,, \; 0<\beta<1\,,\; c>0.$}
\\
{\it Then, with  $\lambda =\rho$  in case} (A),
{\it   $\lambda = c\pi/[\Gamma(\beta+1) \, \sin(\beta \pi)]$    in case} (B)
\\
{\it we have the asymptotics $1-\widetilde \phi(s)\sim \lambda s^\beta$  for $0<s \to 0$.}

	 We will now outline the {\it well-scaled passage to the diffusion limit} by which, 
	  via rescaling space and time in a combined way, we will arrive at the Cauchy problem for 
	  the space-time fractional diffusion equation. 
	  Assuming the conditions of the two lemmata fulfilled. 
	  we carry out this passage in the Fourier-Laplace domain. 
	  For rescaling we multiply the jumps and the waiting times by  positive factors $h$ and $\tau$ 
	   and so obtain a random walk 
	   $x_n(h)=  (X_1 +X_2 +\dots + X_n)\,h $
	     with jump instants $t_n(h)= (T_1+T_2 + \dots + T_n)\, \tau$.  
	   We study this rescaled random walk under the intention to send $h$ and $\tau$  towards 0. 
	   Physically, we change the units of measurement from 1 to $1/h$  in space, from 1 to $1/\tau$  in time, 
	   respectively, making intervals of moderate size numerically small, 
	   and intervals of large size numerically of moderate size, 
	   in this way turning from the microscopic to the macroscopic view. 
	   Noting  the densities 
	   $w_h(x)= w(x/h)/h$ and $\phi _\tau(t/\tau)/\tau$
	   of the reduced jumps and waiting times,      
	   we get the corresponding transforms $\widehat w_h(\kappa)= \widehat (\kappa h)$,
	   $\widetilde \phi _\tau(s)= \widetilde \phi(\tau s) $, 
	   and in analogy to the Montroll-Weiss equation (2.5) the result
$$
\widehat{\widetilde p}_{h,\tau}(\kappa ,s) \!=\! 
\frac{1-\widetilde \phi_\tau (s)}{s} \, \frac{1}{1- \widehat w_h(\kappa )\,\widetilde  \phi _\tau(s)}
\!=\! \frac{1-\widetilde \phi (\tau s)}{s} \, \frac{1}{1- \widehat w(h\kappa )\,\widetilde \phi(\tau s)}\,.
\eqno(3.1) $$
      Fixing now $\kappa$  and $s$ both as $\ne 0$, replacing  
	  $\kappa$ by $h \kappa$  and  $s$ by $\tau s$    in Lemma 1 and Lemma 2, 
	  sending $h$  and $\tau$  to zero, we obtain by a trivial calculation the asymptotics
$$
\widehat{\widetilde p}_{h,\tau}(\kappa ,s) \sim 
\frac{\lambda \tau^\beta s^{\beta-1}}{\mu (h |\kappa|)^\alpha + \lambda (\tau s)^\beta}
\eqno(3.2) $$                   
 that we can rewrite in the form
$$ \widehat{\widetilde p}_{h,\tau}(\kappa ,s) \sim
\frac{ s^{\beta-1}}{r(h,\tau)|\kappa|^\alpha +  s^\beta}
\q \hbox{with}\q r(h,\tau)= \frac{\mu h^\alpha}{\lambda \tau^\beta}\,. \eqno  (3.3)$$          
Choosing  $r(h, \tau)\equiv 1$   (it suffices to choose  $r(h, \tau)\to 1$) we get
$$
\widehat{\widetilde p}_{h,\tau}(\kappa ,s)\to \widehat{\widetilde p}_{0,0}(\kappa ,s)
= \frac{s^{\beta-1}}{|\kappa|^\alpha + s^\beta}\,,
\eqno(3.4) $$
We honour by the name {\it scaling relation} our condition
$$
\frac{\mu h^\alpha}{\lambda \tau^\beta}\equiv 1\,.
\eqno(3.5)$$
Via $\tau= \left(\mu/\lambda)h^\alpha\right)^{1/\beta}$ we can eliminate the parameter $\tau$, 
apply inverse Laplace transform to (3.2), fix $\kappa$ and send $h\to 0$.
 So, by the continuity theorem (for the Fourier transform of a probability distribution, see Feller (1971),                                               .
       we can identify 
	   $$ \widehat{\widetilde p}_{0,0}(\kappa ,s) = 
	 \frac{s^{\beta-1}}{|\kappa|^\alpha + s^\beta}   $$
 as the Fourier-Laplace solution $\widehat{\widetilde u}(\kappa,s)$ of the 
 {\it space-time fractional Cauchy problem} (for  $x\in \RR$, $t\ge 0$) 
$$
{\, _t}D_{*}^{\, \beta }\, u(x,t)
  = 
 {\, _x}D_{0}^{\,\alpha} \,u(x,t)
\,,\; u(x,0)=\delta(x)\,, \; 0<\alpha \le 2\,,\; 0<\beta \le 1\,.
\eqno(3.6)$$              
Here, for $0<\beta \le 1$, we denote by ${\, _t}D_{*}^{\, \beta }$  
the regularized  fractional differential operator, see  Gorenflo and Mainardi (1997), according to
$$   {\, _t}D_{*}^{\, \beta } \, g(t)
={\, _t}D^{\, \beta }\left[g(t)-g(0)\right]
\eqno (3.7)$$           
with the Riemann-Liouville fractional differential operator
$$  
 \,_tD^\beta  \,g(t) :=
 \begin{cases}
 {\ds \rec{\Gamma(1-\beta )}}
  {\ds {d\over dt}\,\int_0^t
    {g(t')\,d\tau  \over (t-t' )^{\beta  }} 
	} \,, \q 0<\beta<1 \,,\\
	  {\ds \frac{d}{dt} g(t)}\,,\q \beta=1\,.
	  \end{cases}
\eqno(3.8)$$
Hence, in longscript:
$$   {\, _t}D_{*}^{\, \beta } \, g(t)
=
\begin{cases}
 {\ds \rec{\Gamma(1-\beta )}}
  {\ds {d\over dt}\,\int_0^t
    {g(t')\,d\tau  \over (t-t' )^{\beta  }} } -
	{\ds \frac{g(0)t^{-\beta}} {\Gamma(1-\beta)}}\,,\q 0<\beta<1 \\
	{\ds \frac{d}{dt} g(t)}\,,\q \beta=1\,.
	\end{cases}
  \eqno(3.9)$$
If $g'(t)$   exists we can write
$${\, _t}D_{*}^{\, \beta } \, g(t) =
{\ds \rec{\Gamma(1-\beta )}}
  {\ds \,\int_0^t
    {g'(t')  \over (t-t' )^{\beta  }}\, dt' } \,, \q 0<\beta<1 \,, $$
	and the regularized fractional derivative coincides with the form introduced by Caputo, see
	Caputo and Mainardi (1971), Gorenflo and Mainardi (1997), Podlubny (1999),
	 henceforth referred to as the Caputo  derivative. 
Observe that in the special case $\beta=1$  
 the two fractional derivatives 
${\ds{\, _t}D_{*}^{\, \beta } \, g(t)}$
   and 
   ${\, _t}D^{\, \beta } \, g(t)$
     coincide, both then being equal to $g'(t)$.

      The Riesz operator   is a pseudo-differential operator according to
$$
\widehat{\,_xD_0^\alpha \, f} =  - |\kappa| ^\alpha   \,
 \widehat f(\kappa) \,,  \q \kappa \in \RR\,,
\eqno(3.10) $$                                       
compare Samko, Kilbas and Marichev (1993) and  Rubin (1996). 
It has the Fourier symbol $-|\kappa|^\alpha$.
      
	    In the transform domain (3.6) means     
$$ s^{\beta-1}\, \widehat{\widetilde u}(\kappa,s) - s^{\beta-1}= -|\kappa|^\alpha\,
 	\widehat{\widetilde u}(\kappa,s)$$ 
		 hence 
$$
\widehat{\widetilde u}(\kappa,s)= \frac{s^\beta-1}{|\kappa|^\alpha + s^\beta}\,, 
\eqno(3.11)$$             
and looking back at (3,4 ) we see:  
${\ds \widehat{\widetilde u}(\kappa,s)} = {\ds \widehat{\widetilde p}_{0,0}(\kappa,s)}$. 
Thus, under the scaling relation  (3.5),  the Fourier-Laplace solution of the CTRW 
integral equation (2.1) converges to the Fourier-Laplace solution of 
the space-time fractional Cauchy problem (3.6), 
and we conclude that the sojourn probability of the CTRW converges weakly (or "in law") 
to the solution of the Cauchy problem for the space-time fractional diffusion equation
for every fixed $t>0$. 
Later in this paper we will present another way of passing to the diffusion limit, 
a way in which by decoupling the transitions in  time and in space we circumvent doubts on the correctness 
of the transition.

      For a comprehensive study of integral representations of the solution to the Cauchy problem (3.6) 
	  we recommend the paper by Mainardi, Luchko and Pagnini (2001).
	  
 \subsection*{Subdiffusive and Superdiffusive Behaviour}
  With regard to the parameters  $\alpha$ and $\beta$  in equation (3.6) we single out 
  the cases (i), (ii) and (iii)  
  by  attributing names to them.  
\\ (i)   $\alpha=1\,,\; \beta=1$:      {\it normal}     or {\it Gaussian diffusion} 
(according to ${ \frac{\d u}{\dt}}= { \frac{\d^2 u}{\dx^2}}$), 
\\ (ii) $\alpha=2\,,\; 0<\beta<1$:    {\it time-fractional diffusion}, 
\\ (iii) $0<\alpha <2\,,\; \beta= 1$:   {\it space-fractional diffusion}.
 
Let us now consider (compare Gorenflo nnd Abdel-Rehim (2004)) the equation (3.6) and its solution (3.11) 
in transform space that describe the evolution of the sojourn probability density
$u(x,t)$  of a wandering particle starting in the origin $x=0$  at the initial instant $t=0$. 
We call this behaviour {\it subdiffusive} if the {\it variance}
$$ \langle (x(t))^2 \rangle \, = \, (\sigma(t))^2 := 
{\ds \int_{-\infty}^{+\infty} \!\! x^2\, u(x,t)\, dx} $$
   behaves for $t\to \infty$ like a power $t^\gamma$   with $0<\gamma <1$ , 
{\it normal} if $\gamma=1$, 
{\it superdiffusive} if $\gamma >1$    or if this variance is infinite for positive $t$. 
Using the fact that  by Fourier transform theory   
$$ {\ds \int_{-\infty}^{+\infty} \!\! x^2\, u(x,t)\, dx} = -
{\ds \frac{\d^2}{\d \kappa^2} \widehat u(\kappa,t)|_{\kappa=0}}\,$$
and writing the right hand side of (3.11) as an infinite  series in powers of   $\kappa/s$
(convergent for $s>1$) we find by termwise Laplace inversion
$$ \widehat u(\kappa,t) = 1- \frac{|\kappa|^\alpha t^\beta}{\Gamma(1+\beta)} 
+ \frac{|\kappa|^{2\alpha} t^{2\beta}}{\Gamma(1 +2\beta)} - +\dots
= E_\beta \left(-|\kappa|^\alpha t^\beta \right)$$
from which for $t>0$ we obtain, for all $0<\beta \le 1$ the result,  
    $$ (\sigma(t))^2 =
	\begin{cases}
	{2 t^\beta}/{\Gamma(1+\beta)}\,, &\q \hbox{if} \q \alpha=2\,,\\
	\infty \,,  &\q \hbox{if} \q 0<\alpha<2\,.
	\end{cases}
	 $$ 
In the special case (i) of Gaussian diffusion this reduces to   
$(\sigma(t))^2 = 2t $.
\section*{4. Thinning (Rarefaction) of a Renewal Process under Power Law Regime}
We  are going to give an account of the essentials of thinning a renewal process 
with power law waiting times, thereby leaning on the presentation by Gnedenko and Kovalenko (1968) 
but for reasons of transparency not decorating the power functions by slowly varying functions. 
	   Compare also Mainardi, Gorenflo and Scalas (2004) and Gorenflo and Mainardi (2008). 
     
Again (as in Section 2) with the $t_n$ in strictly increasing order, the time instants of 
a renewal process, $0=t_0<t_1<t_2<\dots$, with {\it iid} 
waiting times $T_k=t_k-t_{k-1}$  (generically denoted by $T$), 
 {\it thinning} (or {\it rarefaction}) means that for each positive index $k$ a decision is made: 
	 the event happening in the instant $t_k$   
	 is deleted with probability $p$ (where $0<p<1$) 
	 or is maintained with probability $q=1-p$. 
	 This procedure produces a thinned (or rarefied) renewal process, 
	 namely one with fewer events. 
Of particular interest for us is the case of $q$  near zero which results in 
very few events in a moderate span of time. To compensate for this loss 
(wanting to keep a moderate number of events in a moderate span of time) 
we change the unit of time which amounts to multiply the (numerical value of)  
the waiting time with a positive factor $\tau$   
so that we get waiting times $\tau T_k$  and instants $\tau t_k$  
in the rescaled process. 
Loosely speaking, it is our intention to dispose on $\tau$  
in relation to the rarefaction factor $q$ in such way that for very small $q$ 
in some sense the "average" number of events per unit of time remains unchanged. 
We will make these considerations precise in an asymptotic sense.
      
	  Denoting by $F(t) = P(T\le t)$  the probability distribution function 
	  of the original waiting time $T$, 
	  by $f(t)$  its density (generically this density is a generalized function represented 
	  by a measure) so that
	  $F(t) ={\ds \int_0^t f(t')\,dt'}$,
	  and analogously for the functions
	  $F_k(t)$   and $f_k(t)$,  
	  the distribution and density, respectively, of the sum of $k$
	   waiting times, we have recursively
$$ f _1(t)= f(t)\,, \; f_k(t) =\int_0^t f_{k-1}(t)\, dF(t') \quad \hbox{for}\quad k\ge 2\,.
\eqno(4.1)$$          
Observing that after a maintained event of the original process the next one is kept with 
probability $p$ but dropped with probability $q$ in favour of the second-next with probability
$pq$   and,  generally $n-1$  events are dropped  in favour 
of the  $n$-th next with probability $p^{n-1}q$,
we get for the waiting time density of the thinned process the formula
$$ g_q(t)=  \sum_{n=1}^\infty q\, p^{n-1}\, f_n(t)\,.\eqno(4.2)$$                                     
With the modified waiting time $\tau T$   
we have  $P(\tau T \le t)= P(T\le t/\tau)= F(t/\tau)$, 
hence the density  $f(t/\tau)/\tau$, and analogously  for the density of the sum 
of $n$  waiting times $f_n(t/\tau)/\tau$.
 The density of the waiting time of the sum of $n$ waiting times of the 
 rescaled (and thinned) process now turns out as
$$g_{q,\tau}(t)=  \sum_{n=1}^\infty q\, p^{n-1}\, f_n(t/\tau)/\tau\,.\eqno(4.3)$$                      
In the Laplace domain we have $\widetilde f_n(s)= (\widetilde f(s))^n$, hence 
(using $p=1-q$)
 $$ \widetilde g_q(s)=  \sum_{n=1}^\infty q\, p^{n-1}\, (\widetilde f(s))^n =
 \frac{q \, \widetilde f(s)}{1-(1-q)\widetilde f(s)}\,.
 \eqno(4.4)$$                             
By rescaling we get
$$
\widetilde g_{q,\tau}(s)=  \sum_{n=1}^\infty q\, p^{n-1}\, (\widetilde f(\tau s))^n =
 \frac{q \, \widetilde f(\tau s)}{1-(1-q)\widetilde f(\tau s)}\,.
\eqno(4.5)$$                        
Being interested in stronger and stronger thinning ({\it infinite thinning}) 
let us consider a scale of processes  with the parameters  $q$  of  {\it thinning}
 and $\tau$  of {\it rescaling} tending to zero under 
 a {\it scaling relation} $q=q(\tau)$   yet to be specified.
 
      Let us consider two cases for the (original) waiting time distribution, 
	  namely as in Lemma 2 of Section 3 case (A) of a finite mean waiting time and case (B) 
	  of a power law waiting time.  
	  We assume  
$$ \lambda := \int_0^\infty t' \, f(t')\,dt' <\infty \q (A), \q \hbox{setting} \q
\beta=1,\eqno(4.6A)$$
or
$$ \Psi(t)= \int_t^\infty f(t')\, dt' \sim \frac{c}{\beta} t^{-\beta}
\q \hbox{for} \q t \to \infty\q \hbox{with} \q 0<\beta<1\,.  \eqno(4.6B)$$
In case (B) we set   
$$ \lambda = \frac{c \pi}{\Gamma(\beta+1)\, \sin (\beta \pi)}\,. $$
By Lemma 2 of Section 3 we have $\widetilde f(s)= 1-\lambda s^\beta +o(s^\beta)$  
 for   $0<s \to 0$.
  
  Passing now to the limit $q\to 0$  of {\it infinite thinning} under the {\it scaling relation}
$$  q = \lambda \tau ^\beta \eqno(4.7)$$                                                
for {\it fixed} $s$ the Laplace transform  (4.5) of the rescaled density
$g_{q,\tau}(t)$  of the thinned process   tends   to
$\widetilde g(s) =1/(1+s^\beta)$
    corresponding to the Mittag-Leffler density
 $$ g(t)= - \frac{d}{dt} E_\beta(-t^\beta)= \phi^{ML}_\beta(t)\,. \eqno(4.8)$$                                  
  Thus, the thinned process converges weakly to the {\it Mittag-Leffler renewal process}  
  described in  Mainardi, Gorenflo and Scalas (2004) 
  (called {\it fractional Poisson process} in Laskin (2003)) 
  which in the  special case $\beta=1$  
  reduces to the Poisson process. 
  In this sense the Mittag-Leffler renewal process is 
  {\it asymptotically universal for power law renewal processes}.               
 \section*{5. Mittag-Leffler Waiting Time and Space-Time Fractional Diffusion}
     Let us sketch how, under the power law assumptions of Lemma 1 and Lemma 2, 
the {\it simultaneous passage to the limit} in inversion of  the Fourier and Laplace transforms
  in (3.4) can be circumvented.   
  Leaning on our presentations in  Mainardi et al. (2000), Gorenflo et al. (2001)
  and Scalas, Gorenflo and Mainardi (2004), 
  we introduce a {\it memory function} $H(t)$   via which we will arrive 
  at an evolutionary integral equation for the sojourn probability density  $p(x,t)$. 
  For illustration we will soon consider a few special choices 
  for this function. 
  By {\it rescaling} and {\it respeeding} the process in time $t$ 
  and passing to an appropriate limit we will get a time-fractional evolution equation 
  for $p(x,t)$   
  (in fact a time-fractional generalization of the Kolmogorov-Feller equation) 
  that arises also by direct insertion  of  the Mittag-Leffler waiting time density 
  into the CTRW integral equation (2.1) as we can see in Hilfer and Anton (1995).  
  Via a second {\it respeeding}, obtained by  rescaling the spatial variable
  $x$ , we will arrive at the Cauchy problem (3.6) for space-time fractional diffusion. 
  We keep the notations of Sections 2 and 3.
  
        First, we introduce in the Laplace domain the auxiliary function 
$$   \widetilde{H} (s) = \frac{1- \widetilde{\phi}(s) }
{ s\, \widetilde{\phi}(s)}
 = \frac{\widetilde{\Psi}(s) }{\widetilde{\phi}(s)}\,,\eqno(5.1)$$                     
and see by trivial calculation that (2.4) is equivalent to 
$$\widetilde{H} (s) \left[   s \widehat{\widetilde p}(\kappa,s)-1 \right]=
\left[\widehat w(\kappa)-1\right]\widehat{\widetilde p}(\kappa,s)\,,
\eqno(5.2)$$                    
meaning in the space-time domain the {\it generalized Kolmogorov-Feller equation}
   $$  \int_0^t   \!\! H(t-t')\,
 \frac{\d}{\d t'} p(x,t')\, dt'  =   
    -p(x,t)+ \int_{-\infty}^{+\infty}\!\! w(x-x')\, p(x',t)\, dx'\,,
\eqno(5.3) $$
with $p(x,0)=\delta(x)$.

Note that  (5.1) can be inverted to
$$ \widetilde \phi(s)= \frac{1}{1 +s \widetilde H(s)}\,. \eqno(5.4)$$                                        
We may play with equation (5.1), trying special choices for $\widetilde H(s)$  
to obtain via (5.4) meaningful waiting time densities $\phi(t)$. 
Or we chose  $H(t)$ hoping again to get via (5.4) a meaningful density $\phi(t)$.
In accordance with our inclination towards power laws let us take      
	 $\widetilde H(s)= s^{\beta-1}$
	 and distinguish the cases  (i)$\beta=1$    and (ii) $0<\beta<1$. 
\\	 
      Case (i) yields  
	  $$  \widetilde{H} (s)\equiv 1\,,\; H(t)=\delta(t)\,, \;
	  \widetilde \phi(s)= \frac{1}{1+s}\,, \;\phi(t)= \exp (-t)\,,$$
	namely 	  the exponential waiting time density, and (5.2) reduces to 
$$
 s \widehat{\widetilde p}(\kappa,s)-1 =
\left[\widehat w(\kappa)-1\right]\widehat{\widetilde p}(\kappa,s)\,,
\eqno(5.5)$$                                 
in the space-time domain the classical {\it Kolmogorov-Feller equation} 
$$   \frac{\d}{\d t'} p(x,t')\, dt'  =   
    -p(x,t)+ \int_{-\infty}^{+\infty}\!\! w(x-x')\, p(x',t)\, dx'\,,  
	     \eqno(5.6)$$
  with $ p(x,0)=\delta(x)$
  \\  
	   Case (ii) yields   
 $$   H(t)=\frac{t^{-\beta}}{\Gamma(1-\beta)}=: H^{ML}_\beta(t), \;
	  \widetilde \phi(s)= \frac{1}{1+s^\beta},\; \phi(t)= -\frac{d}{dt}E_\beta (-t^\beta)
	  = \phi^{ML}_\beta(t),$$
	  namely 
  the Mittag-Leffler waiting time density introduced in Section 1 by formula (1.4). 
  With the Caputo derivative operator 
  ${\ds \,_tD_*^\beta}$ 
     of (3.7) we get the Cauchy problem
$$    \, _tD_*^\beta \,  p(x,t) =
     -  p(x,t) +   \int_{-\infty}^{+\infty} w(x-x')\,
   p(x',t) \, dx'\,, \eqno(5.7) $$ 
with $ p(x,0)=\delta(x)$.
\\
{\bf Remark:} 
Because of (3.7) the equations (5.7) and (5.6) coincide in the special case $\beta=1$.
\subsection*{Rescaling and Respeeding}
  Let us now manipulate the generalized Kolmogorov-Feller equation (5.3) 
  by working on it in the Laplace domain via (5.2). 
  
  {\it Rescaling} time means: With a positive scaling factor $\tau$ (intended to be small) 
  we replace the waiting time $T$ by $\tau T$. 
  This amounts to replacing the unit 1 of time by $1/\tau$, and if $\tau <<1$  
  then in the rescaled process there will happen very many jumps in a moderate span 
  of time (instead of the original moderate number in a moderate span of time). 
  The rescaled waiting time density and its corresponding Laplace transform
   are $\phi _\tau(t)= \phi(t/\tau)/\tau$, $ \widetilde \phi _\tau(s)=\widetilde \phi(\tau s) $.   
  Furthermore:
  $$  \widetilde{H}_\tau (s)= \frac{1-\widetilde \phi _\tau(s)}{s\, \widetilde \phi _\tau (s) }=
   \frac{1- \widetilde{\phi}(\tau s)}{s\, \widetilde{\phi}(\tau s)}
    \,,\; \hbox{hence}\;
\widetilde \phi _\tau(s)
= \frac{1}{1 +s\, \widetilde H_\tau(s)}\,,
\eqno(5.8)$$     
and (5.2) goes over into 
$$
\widetilde{H}_\tau  (s) \, \left[
  s\widehat{\widetilde {p}}_{\tau}(\kappa ,s)-1\right] =
 \left[ \widehat w(\kappa) - 1\right]\,
   \widehat{\widetilde {p}}_{\tau} (\kappa ,s)\,.
\eqno (5.9)$$  
\\              
 {\bf Remark:} Note that in this Section 5 the position and meaning of the indices attached 
 to the generic density $p$ are convenient but {\it different} from those in Section 3.

   {\it Respeeding} the process means multiplying the left hand side 
   (actually ${\ds \frac{\d}{\d t'}p(x,t')}$  of equation (5.3) by a positive factor $1/a$, 
   or equivalently its right hand side by a positive factor $a$. 
   We honour the number $a$ by the name {\it respeeding factor}. 
   $a>1$   means {\it acceleration},  $a<1$ {\it deceleration}. 
   In the Fourier-Laplace domain the rescaled and respeeded CTRW process then assumes the form, 
   analogous to  (5.2) and (5.9),
$$\widetilde{H}_{\tau,a} (s) \left[   s {\widehat{\widetilde p}}_{\tau,a}(\kappa,s)-1 \right]
=
a\, \left[\widehat w(\kappa)-1\right]{\widehat{\widetilde p}}_{\tau,a}(\kappa,s)\,,
\eqno(5.10)$$                
  with    
  $$ \widetilde{H}_{\tau,a} (s) = \frac{\widetilde{H}_{\tau} (s)}{a}=
  \frac{1-\widetilde\phi(\tau s)}{a\,s\,\widetilde \phi(\tau s)}\,,$$ 
What is the effect of such {\it combined rescaling and respeeding}? We find 
$$\widetilde \phi _{\tau,a}(s)= \frac{1}{1+ s \widetilde H _{\tau,a}(s)}=
\frac{a\, \widetilde\phi(\tau s)}{1- (1-a)\,\widetilde \phi(\tau s)}\,,
\eqno (5.11)$$                          
and are now in the position to address the
{\bf Asymptotic universality of the Mittag-Leffler waiting time density}.
\\ 
       Using Lemma 2 with  $\tau s$  in place of $s$ and taking 
$$ a = \lambda \tau^\beta\,, \eqno(5.12)$$                              
fixing $s$ as required by the continuity theorem of probability for Laplace transforms, 
the asymptotics  
$\widetilde \phi(\tau s) = 1 -\lambda (\tau s)^\beta) +o((\tau s)^\beta)$  for $\tau \to 0$  implies 
$$
\widetilde\phi _{\tau, \lambda \tau^\beta} (s )\!=\!
\frac {\lambda \tau^\beta \left[ 1-  \lambda (\tau s)^\beta) +o((\tau s)^\beta)\right]}
{ 1 - (1-\lambda \tau^\beta)\left[ 1-  \lambda (\tau s)^\beta) +o((\tau s)^\beta)\right]  }
\to \frac{1}{1+s^\beta}\!=\! \widetilde \phi^{ML}_\beta\,,
\eqno(5.13)$$
corresponding to the Mittag-Leffler density   
$${\ds \phi^{ML}_\beta(t)} = {\ds - \frac{d}{dt}E_\beta(-t^\beta)}\,.$$

        Observe that the parameter  $\lambda$ does not appear in the limit
		$1/(1+s^\beta)$. 
		We can make it reappear by choosing the respeeding factor $\tau^\beta$ in place 
		of $\lambda \tau^\beta$. 
		In fact: 
		$$ \widetilde \phi_{\tau, \tau^\beta}   \to \frac{1}{1+\lambda s^\beta}\,. $$             
  Formula (5.13) says that the general density $\phi(t)$  with power law asymptotics 
	  as in Lemma 2 is gradually deformed into the Mittag-Leffler waiting time density
	  $\phi^{ML}_\beta(t)$. 
	  It means that with larger and larger unit of time (by sending $\tau \to 0$) 
	  and stronger and stronger deceleration (by $a=\lambda \tau^\beta$) 
	  as described our process becomes indistinguishable from one 
	  with Mittag-Leffler waiting time (the probability 
	  distribution of jumps remaining unchanged). 
	  Likewise a pure renewal process with asymptotic 
	  power law density becomes indistinguishable 
	  from the one with Mittag-Leffler waiting time 
	  (the fractional generalization of the Poisson process by Laskin (2003) and 
	   Mainardi, Gorenflo and Scalas (2004)). 
	  In fact, we can consider the pure renewal process as a CTRW with 
	  jump density $w(x)=\delta(x-1)$, 
	  the position of the wandering particle representing the counting number 
	  (the number of events up to and including the instant $t$).
\\	  
   {\bf  Remark:} It is instructive to look at the effect of 
   {\it combined rescaling and respeeding} on  the  Mittag-Leffler density    
   $\phi_\beta^{ML}(t)$ itself which by (1.6) also obeys the asymptotic  conditions of Lemma 2. 
   We have
   $\widetilde \phi_\beta^{ML}(s) = 1-s^\beta + o(s^\beta)$
     for $s\to 0$, and  with (5.11) we find the relation
$$ \left( \widetilde \phi_\beta^{ML}\right)_{\tau,a}(s) = 
\widetilde \phi_\beta^{ML}\left(\tau s/a^{1/\beta}\right)\q \hbox{for all} \q \tau>0\,,\; a>0\,,
\eqno(5.14)$$                                   
 expressing the self-similarity of the Mittag-Leffler density. 
 In particular we have the  formulas of self-similarity and invariance
     $$  \left(\phi_\beta^{ML}\right)_{\tau, \tau^\beta}(t)= \phi _\beta^{ML}(t)
	 \q \hbox{for all} \q \tau>0\,,
	 \eqno (5.15)$$         
telling us that the Mittag-Leffler density   is invariant under the transformation (5.11) 
with the respeeding factor $a=\tau^\beta$  in place of $a = \lambda \, \tau^\beta$.

\subsection*{Diffusion Limit in Space}
       In addition to rescaling time we now rescale also the spatial 
	   variable $x$, by replacing the jumps $X$ by jumps $hX$ 
	    with positive scaling factor $h$, intended to be small. 
		The rescaled jump density turns out as $w_h(x)= w(x/h)/h$, 
		corresponding to $\widehat w_h(\kappa)= \widehat w(\kappa h)$. 
		Starting from the Fourier-Laplace representation (5.2) of our CTRW 
		with general waiting time density, 
		we accelerate the spatially rescaled process by the respeeding
		 factor $1/(\mu h^\alpha)$   with $\mu >0$  
		and arrive (using  $q_h$ as new dependent variable) at the equation
$$
\widetilde{H} (s) \, \left[
  s\widehat{\widetilde {q}}_{h}(\kappa ,s)-1\right] =
\frac{\widetilde w(\kappa h)-1}{\mu h^\alpha}\,
   \widehat{\widetilde {q}}_{h} (\kappa ,s)\,,
\eqno(5.16)$$                     
Assuming now the power law assumptions of Lemma 1 satisfied, 
fixing $\kappa$  as required by the continuity theorem for the Fourier transform, 
we get $(\widetilde w(\kappa h)-1)/(\mu h^\alpha) \to -|\kappa|^\alpha$
    for $h\to 0$, 
and writing $u$ in place of $q_0$, in the limit  
$$
\widetilde{H} (s) \, \left[
  s\widehat{\widetilde {u}}(\kappa ,s)-1\right] = - |\kappa|^\alpha
  \,\widehat{\widetilde {u}}(\kappa ,s)\,.
\eqno(5.17)$$      
By Fourier inversion we get   
$$\widetilde{H} (s) \, \left[
  s \widetilde {u}(x ,s)-\delta(x)\right]
  = {\, _x}D_{0}^{\,\alpha} \,\widetilde u(x,s)\,, $$
 and then by Laplace inversion in the space-time 
domain the Cauchy problem
$$
\int_0^t \!\!H (t-t') \, \frac {\d}{\d t} u(x ,t)
  = {\, _x}D_{0}^{\,\alpha} \, u(x,t)\,, \; u(x,0)=\delta(x)\,, \; 0<\alpha \le 2\,.
\eqno(5.18) $$
As in Section 3  we mean by $\,_xD_0^\alpha$ the Riesz pseudo-differential operator 
with Fourier symbol $-|\kappa|^\alpha$  according to formula (3.10).

 Finally, inserting into (5.18) the Mittag-Leffler memory function
    $$ H^{ML}_\beta(t) = 
	\begin{cases}
	{\ds \frac{t^{-\beta}}{\Gamma(1-\beta)}}\,;\hbox{if}\; 0<\beta<1,\\
	 \delta(t)\,; \hbox{if}\; \beta=1 ,
	 \end{cases}
	 $$
we recover the Cauchy problem (3.6) for the space-time fractional diffusion equation,
 namely
$$     {\, _t}D_{*}^{\, \beta }\, u(x,t)
 \, = \, {\, _x}D_{0}^{\,\alpha} \,u(x,t)\,, 
 \; u(x,0)=\delta(x)\,,
  \; 0<\alpha \le 2\,, \; 0<\beta \leq 1
\,.(5.19)$$          
\\
{\bf Comments:} In this Section 5 we have split the passage 
to the limit into a temporal one $\tau \to 0$    and a spatial one
$h\to 0$. 
In Section 3 by the well-scaled (combined) passage to the limit 
we have avoided the concept of respeeding but have 
had to eliminate the parameter $\tau$ 
via the scaling relation (3.5).
But where has the scaling relation (3.5) gone? 
We rediscover it as hidden in the deceleration (5.10) with (5.12)  
and the compensating acceleration (5.16). 
We recommend to compare this splitting with the more abstract technique 
of triangular arrays applied by Meerschaert and Scheffler (2004) and  (2008). 
We think that our method offers intuitive insight into 
the meaning of passing to the diffusion limit.      
\section*{6. Time-Fractional Drift and Subordination}
       In Section 2 we have seen in which way a CTRW is 
	   subordinated to a renewal process , and in Section 5 we have worked out 
	   the effect of subordination under the Mittag-Leffler renewal process, 
	   see equation (5.7). 
	   For the following considerations we hint to the references 
	   Mainardi, Gorenflo and Scalas (2004), Meerschaert et al (2002),
	   Mainardi, Pagnini and Gorenflo (2003), Mainardi, Luchko and Pagnini (2001),
	   Gorenflo, Mainardi and Vivoli (2007), Gorenflo and Mainardi (2008),
	   and again to the papers by Meerschaert and Scheffler (2004) and  (2008).
	   
	   Our aim is to pass to a meaningful limit of the Mittag-Leffler renewal process
	   whose waiting density  is $\phi _\beta^{ML}(t)$ under the restriction $0<\beta <1$,
	   meaning exclusion of the limiting case $\beta=1$ of exponential waiting time.
       
	     Keeping our earlier notations, in particular $T_k$  for the waiting times,  
		$t_k$ for the jump instants, and denoting by
		$N=N(t)= \hbox{max}\, \{k|t_k \le t\}$  the number of renewal or 
		 jump events up to instant $t$ (the "counting number") 
		 we have for a general renewal process the probability
$$ 
P (N(t)=k ) = P(t_k \le t\,, \, t_{k+1}>t)\,.
\eqno(6.1)$$                         
Finding it convenient to embed the formalism of renewal into the CTRW formalism 
we introduce a pseudo-spatial variable  $r$  taking the  values of the counting number 
$N$ and denote the sojourn probability by $q(r,t)$  (with $r>0$, $t\ge 0$), 
by $v(r)$ the jump density. 
Taking (because $N$ runs through the non-negative integers) for the  jump-width 
of the "random walk" so generated the constant 1 we have
$v(r)=\delta(r-1)$   and $\widehat v(\kappa)= \exp (i\kappa)$, and (5.2) yields
$$
\widetilde{H} (s) \, \left[
  s\widehat{\widetilde {q}}(\kappa ,s)-1\right] = \left[\widehat v(\kappa) -1 \right]
  \,\widehat{\widetilde {q}}(\kappa ,s)\,.
\eqno(6.2)$$                        
With the  waiting time density $\phi^{ML}_\beta(t)$
  and correspondingly  
  $$ {\widetilde \phi}^{ML}_\beta(s) = \frac{1}{1+s^\beta}\,, \q
  \widetilde H(s) = s^{\beta-1}\, $$   
  (6.2) goes over in
$$ s^{\beta-1}\, \left[
  s\widehat{\widetilde {q}}(\kappa ,s)-1\right] =
   \left[\exp (i\kappa) -1 \right]  \,\widehat{\widetilde {q}}(\kappa ,s)\,.
\eqno(6.3)$$                       
As we have done in Section 5 for the general CTRW we now ask what 
happens when we pass to the diffusion limit in "space" for the Mittag-Leffler 
renewal process in the CTRW formalism. 
Multiplying the jumps by a positive scaling factor  $\delta$, 
decorating $q$ by such index, replacing $\widetilde v(\kappa)$  
by $\widetilde v(\kappa \delta)$   according to
$ v_\delta(r)= v(r/\delta)/\delta$,  fixing  $\kappa$, 
finally accelerating by applying the factor $\delta^{-1}$  
to the right hand side we obtain the equation
$$ s^{\beta-1}\, \left[
  s\widehat{\widetilde {q}}_\delta(\kappa ,s)-1\right] =
   \delta^{-1}\,\left[\exp (i\kappa \delta) -1 \right]  \,\widehat{\widetilde {q}}_\delta(\kappa ,s)\,. $$                          $$  $$
and  $\delta \to 0$ yields the equation
$$
s^{\beta-1}\, \left[
  s\widehat{\widetilde {q}}_0(\kappa ,s)-1\right] =
   i \kappa \,\widehat{\widetilde {q}}_0(\kappa ,s) 
\eqno (6.4)$$                               
which implies
$$\widehat{\widetilde {q}}_0(\kappa ,s)=\frac{s^{\beta-1}}{s^\beta -i\kappa}\,.
\eqno (6.5)$$                                       
We remark that an analogous limit can by proper scaling  
also be obtained directly from the generic power law renewal process that 
we have discussed in Section 4 on thinning.
 
       Note that (6.4) corresponds to the Cauchy problem for the 
	   positively oriented {\it fractional drift} equation with the Caputo
	   derivative operator  $\, _tD_*^\beta$,
$$
\, _tD_*^\beta q_0(r,t)= -\frac{\d}{\d r}q_0(r,t)\,.
\eqno(6.6)$$     
      Without inverting the transforms in (6.5) we can 
	  recognize the {\it self-similarity} of the function $q_0(r,t)$.
	  With any positive constants $a$ and $b$,  generic functions
	   $f$ and $g$ 
	   we have the correspondence of $f(ax)$   to $\widehat f(\kappa/a)/a$ 
	    and the correspondence of $g(bt)$  to $\widetilde g(s/b)/b$, 
	   hence via (6.5)
	   $ q_0(ax, bt) =  b^{-\beta}\, q_0(ax/b^\beta, t)$,
	    and with $Q_0(x)= q_0(x,1)$  and the similarity variable $x/t^\beta$ 
		 we obtain 
$$ q_0(x,t) = t^{-\beta}\, Q_0(x/t^\beta)\,.
\eqno(6.7)$$                                                
Fourier inversion of (6.5 ) gives 
$$\widetilde q_0(r,s)=
\begin{cases}
{\ds s^{\beta -1}\, \exp(-rs^\beta)\,,} \; \hbox{for}\; r>0\, \\
0\,,\; \hbox{for}\; r<0
\end{cases}
\eqno(6.8)$$
Using the fact that $\exp(-s^\beta)$  is the Laplace transform of the extreme positive-oriented unilateral 
stable density of order $\beta$, namely of $L_\beta^{-\beta}(t)$, we get 
$$ q_0(r,t) = r^{-1/\beta}\, _tJ^{1-\beta}\, L_\beta^{-\beta}\left(t r^{-1/\beta}\right)\,,
 \eqno (6.9)$$                                  
with the Riemann-Liouville fractional integration  
$$\, _tJ^{\gamma} \, g(t)= \frac{1}{\Gamma(\gamma)} \, 
\int_0^t \!\! (t-t')^{\gamma- 1} \, g(t')\, dt'\,,
\q \gamma >0 \,. $$
This solution can be expressed in alternative ways for which we refer to the  references 
cited at begin of this Section. 
We have, for $t>0$,
$$ q_0(r,t) = \frac{t}{\beta}\, r^{-1-1/\beta}\,L_\beta^{-\beta}(t r ^{-1/\beta}) =
  t^{-\beta} \, M_\beta(r t^{-\beta})\,,
\eqno(6.10)$$              
with the $M$-Wright function
$$ M_\beta  (z) \!=\!
  \sum_{n=0}^{\infty}\,
  {(-z)^n\over n!\,\Gamma[-\beta  n + (1-\beta )]} 
  \!=\!
  \rec{\pi}\, \sum_{n=1}^{\infty}\,{(-z)^{n-1} \over (n-1)!}\,
  \Gamma(\beta  n)  \,\sin (\pi \beta   n).$$
  for whose properties and use we refer to  Mainardi, Mura and Pagnini (2010).
  
      We recognize the function $q_0(r,t)$ as the sojourn probability density 
	  of the directing process (the {\it subordinator}) for space-time fractional diffusion (3.6), 
	  according to Meerschaert et al. (2002), Mainardi, Pagnini and Gorenflo (2003),
	  Gorenflo, Mainardi and Vivoli (2007).  
	  
	  Let us use the notations of Gorenflo, Mainardi and Vivoli (2007). 
	  There it is outlined how the random position $x(t)$   of 
	  the wandering particle can be expressed via 
	  an {\it operational time} $t_*$  (the properly scaled limit of the counting number of the renewal process 
	  to which he CTRW is subordinated) from  which the physical time 
	  $t =t(t_*)$   is generated  
	  by a positive-oriented stable extremal process of order $\beta$  whereas 
	  the spatial position is generated by a stochastic process 
	 $y=y(t_*)$, stable of order $\alpha$. 
	  Then we have the process $x=y(t_*(t))$ in {\it parametric representation} 
 $$ t=t(t_*)\,, \q x =y(t_*)\,.$$
In practice of simulation this representation can be used to produce "flash-light shots" 
of a set of points in the $(t,x)$  plane  showing a sequence of positions of a wandering particle. 
One only needs routines for generating random numbers from the relevant stable probability laws. 
As a consequence of infinite divisibility and Markovianity of these laws the sequence of points 
so produced constitutes a sequence of true particle positions 
(possible positions of a really happening process), see Gorenflo, Mainardi and Vivoli (2007).
  
      However, it is more common, see Feller (1971), to treat {\it subordination} directly 
	  in the form $x=x(t)=y(t_*(t))$, 
	  namely to produce first the process
	  $t_*=t_*(t)$  of generating the operational time $t_*$     
	  from the physical time $t$. 
	  
	  The processes $t=t_*(t)$  and $t_*=t_*(t)$
	   are inverse to each other, but if $0<\beta<1$  
	   the   latter 
	  is neither Markovian nor infinitely divisible. 
	  Its probability density function is $q_0(t_*,t)$, evolving in physical time $t$ 
	  and  given by formula (6.10). 
      Sometimes the process  is called a {\it Mittag-Leffler process} 
	  (see Meerschaert et al. (2002) and search for further details also in several chapters 
	  of Feller's famous book). 
	  The motivation for calling the process   after Mittag-Leffler lies in its 
	  various relation to the Mittag-Leffler function. 
	  The inverse Laplace transform of (6.5) is  
	  $\widehat q_0(\kappa,t)= E_\beta(i\kappa t)$ 
	  (a Mittag-Leffler function with imaginary argument) and 
	  the $M$-Wright $M_\beta$  has by Fourier and Laplace transform 
	  two more connections to the Mittag-Leffler function, namely, 
	  see Mainardi, Mura and Pagnini (2010). 
	  For $0<\beta <1$  we have 
	  $\widetilde M_\beta(s)= E_\beta(-s)$
	    and the Fourier correspondence of 
		$M_\beta(|x|)$  to $2\, E_{2\beta}(-\kappa^2)$. 
	  Furthermore, see 
	   Bondesson, Kristiansen and Steutel (1996) and  Meerschaert and Scheffler (2004),
	   the Laplace-type integral
      $$ \int_0^\infty \!\! \exp(-yr)\, q_0(r,t)\, dr = E_\beta(-yt^\beta)\,.$$                                                     
From (6.5) follows by Fourier inversion 
$$ q_0(r,t)= \frac{1}{2\pi}\, 
VP \,\int_{-\infty}^{+\infty} \!\! \exp(-i\kappa x)\, E_\beta(i\kappa t^\beta)\, d\kappa\,.
\eqno(6.11)$$                     
In the case $0<\beta<1$   the probabilty law of the process $t_*=t_*(t)$ is for no positive $t$ 
infinitely divisible, see Bondesson, Kristiansen and Steutel (1996) and  Meerschaert and Scheffler (2004).
This process must not be confused with the process 
that is called  {\it Mittag-Leffler process} by Pillai (1990) 
and that is obtained via the infinite divisibility of the Mittag-Leffler distribution 
whose density is 
$\phi _\beta^{ML} = - \frac{d}{dt} E_\beta(-t^\beta)$. 

We display the resulting subordination formula 
(compare Meerschaert et al. (2002) and Gorenflo, Mainardi and Vivoli (2007)) 
for the solution $u(x,t)$  of the Cauchy problem (3.6) (in formula (5.19) repeated):        
$$
u(x,t)= \int_0^\infty \!\! f_\alpha(x,r)\, q_0(r,t)\, dr
\eqno (6.12)$$                                        
where the symmetric stable density  $f_\alpha(x,r)$ has the Fourier transform   
$\widehat f_\alpha (\kappa, r) \,\exp(-r |\kappa|^\alpha)$
and the variable $r$  represents the operational time $t_*$.

Let us finally observe that in the excluded limiting case $\beta=1$  by again identifying 
$r$ with  $t_*$,  equation (6.8)  leads  to  the result
$ q_0(t_*,t)= \delta(t-t_*)$, hence $t=t_*$ which
         means  that in the case of exponential waiting time the physical time and the 
         operational time coincide.  
\section*{7. Conclusions}
     We have discussed some (essentially two) ways of passing to the diffusion 
	 limit from continuous time random walk with power laws in time and in space 
	 (for transparency of presentation spatially symmetric and one-dimensional), 
	 namely
	 \\
	  (i) what we call {\it well-scaled passage to the limit} where the rescalings of 
	  time and space are carried out in a combined way, 
	\\
	(ii) carrying out the passages {\it separately in time and in space}.
	 
	The limit in time, by aid of a convenient memory function, 
	leads to the Mittag-Leffler waiting time renewal process or {\it fractional Poisson process}, 
	and the {\it Mittag-Leffler function} becomes essential for description of long-time behaviour 
	(of renewal processes and of CTRW) whereas 
	for the {\it wide-space view} stable distributions take the role. 
	In the time domain there are two passages to the limit. 
	The first one leads to an extremal stable density evolving in time, 
	the other one by condensing the corresponding counting process to smaller and smaller 
	counting-step-size leading to a Mittag-Leffler process as the 
	subordinator of the continuous time random walk.
	 
	Actually we have obtained this subordinator by another splitting: by passing first 
	to the Fractional Poisson process and from this then to the subordinator. 
	But this additional splitting can be avoided. 
	For our way of analyzing the transition to the limit in power law renewal processes 
	we have got inspiration from studying the theory of thinning  
	such processes, and we have discovered the important analogy of a limit formula 
	in the Laplace domain. 
	
	Quite generally, for performing the necessary investigations in CTRW theory 
	the Fourier-Laplace domain is the most convenient operational playground since 
	a long time. 
	To make well visible the basic ideas we have avoided measure-theoretic 
	and functional-analytic terminology, 
	hoping so to be not too difficult for people not so well trained in these fields.
 
 \section*{Acknowledgement} 
 This lecture is based on joint research with Professor F. Mainardi, 
 to whom the author is grateful for long-lasting collaboration.
 
\section*{References}                                           
Abramowitz, M. and I. Stegun, I. (1965): 
{\it Handbook of Mathematical Functions}, Dover Publications, New York. 
\vsn 
Balakrishnan, V. (1985): 
Anomalous diffusion in one dimension, {\it Physica A} {\bf 32}, 569--580.
\vsn 
Bingham, N.H,  Goldie, C.M  and J Teugels, J.L (1967): 
{\it Regular Variation}, Cambridge University Press, Cambridge.  
\vsn 
Bondesson, L., Kristiansen, G.K. and  Steutel, F.W. (1996): 
Infinite divisibility of random variables and their integer parts, 
{\it Statistics and Probability Letters}  {\bf 28},  271--278.
\vsn 
Caputo, M. and  Mainardi, F (1971): 
Linear models of dissipation in an anelastic solid, 
{\it Riv. Nuovo Cimento (Ser. II)} {\bf 1},  161--198. 
\vsn 
Cox, D.R. (1967): 
{\it Renewal Theory}, Second edition, Methuen. London. 
\vsn 
 Erd\'elyi, A., Magnus, W.,  Oberhettinger, F. and  Tricomi, F.G. (1955): 
{\it Higher Transcendental Functions},  Bateman Project, McGraw-Hill, New York, 
Vol. 3, Chapter. 18: Miscellaneous Functions, pp. 206--227.
\vsn 
 Feller, W. (1971): {\it An Introduction to Probability Theory and its Applications}, Vol. II.  
Second Edition. Wiley, New York. 
\vsn 
Gel'fand and  Shilov, G.E. (1964): 
{\it Generalized Functions}, Volume I, Properties and Operations. 
Academic Press, New York and London. 
\vsn 
Gnedenko, B.V. and  Kolmogorov, A.N. (1954): 
{\it Limit Distributions for Sums of Independent Random Variables}, 
Addison-Wesley, Cambridge/Mass..  Translated from Russian.
\vsn 
Gnedenko B.V. and  Kovalenko, I.N. (1968): 
{\it Introduction to Queueing Theory}, 
Israel Program for Scientific Translations, Jerusalem. Translated from Russian. 
\vsn 
Gorenflo, R. and  Abdel-Rehim, E.A. (2004): 
From power laws to fractional diffusion: the direct way, 
{\it Vietnam Journal of Mathematics} {\bf 32 (SI)},   65--75.
\vsn 
Gorenflo, R.,  Loutschko, J.  and   Luchko, Y. (2002): 
Computation of the Mittag--Leffler function     and its derivative,  
{\it Fractional Calculus and Applied Analysis} {\bf  5},  491--518. 
\vsn 
Gorenflo, R. and Mainardi, F. (1997): 
Fractional Calculus: Integral and Differential Equations of Fractional Order,  
In: A. Carpinteri and F. Mainardi (Editors): 
Fractals and Fractional Calculus in Continuum Mechanics. Springer, Wien and New York,
 pp. 223-276.
 E-print in http://arxiv.org/abs/0805.3823
\vsn 
Gorenflo, R. and  Mainardi, F. (2008): 
Continuous time random walk, Mittag-Leffler waiting  time and fractional diffusion: 
mathematical aspects,   in: 
 Klages, R,  Radons, G., and Sokolov, I.M. (Editors):  
{\it Anomalous  Transport, Foundations and Applications}, 
 Wiley-VCH Verlag Gmbh \& Co.  KGaA, Weinheim, Germany,  pp. 93--127. 
 E-print arXiv:cond-mat/07050797.
 \vsn   
 Gorenflo, R. and  Mainardi, F. (2009): 
Some recent advances in theory and simulation of fractional diffusion processes,  
{\it Journal of Computational  and Applied Mathematics}  {\bf 229}, 400--415.
\vsn 
Gorenflo, R., Mainardi, F., Scalas, E. and M Raberto, M.(2001): 
Fractional calculus and continuous-time finance III: the diffusion limit. 
In:  Kohlmann, M. and  Tang, S. (Editors): 
{\it Mathematical Finance}, Birkhäuser: Basel, Boston, Berlin, pp 171--180.
\vsn 
Gorenflo, R., Mainardi, F. and  Vivoli, A. (2007):  
Continuous-time random walk and parametric subordination in fractional diffusion, 
{\it Chaos, Solitons and Fractals} {\bf 34}, 87--103.
\vsn  
Gorenflo, R. and  Vivoli, A. (2003): 
Fully discrete random walks for space-time fractional diffusion equations, 
{\it Signal Processing} {\bf 83},  2411--2420.  
 \vsn 
 Gradshteyn, L.S. and  Ryzhik, I.M. (2000): 
{\it Tables of Integrals, Series and Products}, Sixth Edition,  
Academic Press, San Diego. Translated from the Russian.
\vsn 
Kilbas, A.A.,  Srivastava, H.M. and  Trujillo, J.J. (2006): 
{\it Theory and Applications of Fractional Differential Equations}, Elsevier, Amsterdam. 
\vsn 
Haubold, H.J., Mathai, A.M. and Saxena, R.K.(2009): 
Mittag-Leffler functions and their applications. 
E-print arXiv:0909.0230v2. 
 \vsn 
 Hilfer, H. and  Anton, L. (1995): 
Fractional master equations and fractal time random walks, 
{\it Physical Review E} {\bf 51},  R848--R851.
\vsn 
Hille, E. and  Tamarkin, J.D. (1930): 
On the theory of linear integral equations, 
{\it Annals of Mathematics} {\bf 31},  479--528.
\vsn 
Laskin, N. (2003): 
Fractional Poisson processes, 
{\it Communications in Nonlinear Science and Numerical Simulation} {\bf 8},  201--213. 
\vsn 
 Mainardi, F., Gorenflo, R. and  Scalas, E. (2004): 
A fractional generalization of  the Poisson process, 
{\it Vietnam Journal of Mathematics} {\bf 32 (SI)}, 53--64.
  E-print://arxiv.org/abs/math/0701454
\vsn 
Mainardi, F., Luchko, Yu. and Pagnini, G. (2001): 
The fundamental solution of the space-time fractional diffusion equation, 
{\it Fractional Calculus and Applied Analysis} {\bf 4},  153--192.
\vsn 
Mainardi, F.,  Mura, A. and Pagnini, G. (2010): 
The M-Wright function in time-fractional diffusion processes: a tutorial survey, 
{\it International Journal of Differential Equations},
Vol. 2010, Article ID 104505, 29 pages, 
  http://www.hindawi.com/a104505.html
Electronic Journal published by Hindawi Publishing Corporation, 
 see  http://www.hindawi.com/journals/ijde/contents.html
\vsn 
Mainardi, F.,  Pagnini, G. and  Gorenflo, R. (2003): 
Mellin transform and subordination laws in fractional diffusion processes, 
{\it Fractional Calculus and Applied Analysis} {\bf 6},  441--459.
\vsn 
Mainardi, F., Raberto, M.,   Gorenflo, R.  and Scalas, R. (2000):  
Fractional Calculus and continuous-time finance II: the waiting-time distribution,  
{\it Physica A}  {\bf 287},  468--481.
\vsn 
Mathai, A.M. and  Haubold, H.J. (2008): 
{\it Special Functions for Applied Scientists}, Springer, New York. 
\vsn 
Meerschaert, M.M., Benson, D.A.,  Scheffler, H.-P. and  Baeumer, B. (2002): 
Stochastic solution of space-time fractional diffusion equation, 
{\it Physical Review E} {\bf 65},  041103 (R), 1--4.
\vsn 
Meerschaert, M.M. and Scheffler, H.-P. (2001): 
{\it Limit Theorems for Sums of Independent Random Variables, Heavy Tails in Theory and Practice}, 
Wiley, New York. 
\vsn 
Meerschaert, M.M. and Scheffler, H.-P. (2004):  
Limit theorems for continuous-time random walks with infinite mean waiting times. 
{\it J. Appl. Prob.} {\bf 41},  623--638.
\vsn 
Meerschaert, M.M. and Scheffler, H.-P. (2008):
Triangular array limits for continuous time random walks, 
{\it Stochastic Processes and Their Applications} {\bf 118},  1606--1633.
\vsn 
 Metzler, R. and  Nonnenmacher, Th. F. (2002):  
Space and time fractional diffusion and wave equations, 
fractional Fokker-Planck equations, and physical motivation, 
{\it Chemical Physics} {\bf  284},  67--90.
\vsn 
 Miller, K.S. and  Ross, B. (1993): 
{\it An Introduction to the Fractional Calculus and Fractional Differential Equations}, 
Wiley, New York. 
\vsn 
Mittag-Leffler, G.M. (1903): 
Sur la nouvelle fonction $E_\alpha(x)$, {\it C. R. Acad. Sci. Paris (Ser. II)} 
{\bf 137},   537--539.
\vsn 
Montroll, W.W. and Scher, H. (1973): 
Random walks on lattices, IV: Continuous-time walks and influence of absorbing boundaries, 
{\it J. Stat. Phys.} {\bf 9},  101--135.
\vsn 
Montroll, W.W. and  Weiss, G.H. (1965): 
Random walks on lattices , II,  {\it J. Math. Phys.} {\bf  6}, 167--181.
\vsn  
Newman, M.E.J. (2005): 
Power laws, Pareto distributions and Zipf's law, 
{\it Contemporary Physics} {\bf 46 (5)}, 323--351. 
E-print arXiv:cond-mat/0412004 v3,   29  May 2006. 
\vsn 
 Nonnenmacher, Th.F. (1991): 
Fractional relaxation equations for viscoelasticity and related phenomena,
in:  Casaz-Vazques, J. and Jou, D. (Editors),
{\it Rheological Modelling: Thermodynamical and Statistical Approaches},
Springer Verlag, Berlin, pp. 309--320.
[Lecture Notes in Physics 381] 
\vsn  
Pillai, R.N. (1990): 
On Mittag-Leffler functions and related distributions,
 {\it Ann. Inst. Statist. Math.} {\bf 42},  157--161. 
\vsn 
Podlubny, I (1999): 
{\it Fractional Differential  Equations}, 
Academic Press, San Diego.
\vsn 
Rubin, B. (1996): 
{\it Fractional Integrals and Potentials} 
 Addison Wesley and  Longman, Harlow.
 [Pitman Monographs and Surveys in Pure and Applied Mathematics 82] 
\vsn 
Samko, S.G.,  Kilbas, A.A.,  and  Marichev, G.I (1993): 
{\it Fractional Integrals and Derivatives, Theory and Applications}, Gordon  and Breach, 
New York. Translated from the Russian Edition (Minsk 1987).
\vsn 
 Scalas, E.,   Gorenflo, R.  and  Mainardi, F. (2004):
    Uncoupled continuous-time random walks:
Solution and limiting behavior of the master equation,
{\it Phys. Rev. E} {\bf 69},  011107--1/8.
\vsn 
Seybold, H.  and  Hilfer, R. (2008): 
 Numerical algorithm for calculating the generalized Mittag-Leffler function, 
 {\it SIAM  J. Numer. Anal.} {\bf 47}, 69--88.
\vsn 
Weiss, G.H. (1994): 
{\it Aspects and Applications of the Random Walk}, North-Holland, Amsterdam. 
\vsn 
Wiman, A. (1905): 
(a) \"Uber den Fundamentalsatz in der Theorie der Funktionen,
{\it Acta Mathematica} {\bf 29}, 191--201. 
(b) \"Uber die Nullstellen der Funktionen, {\it Acta Mathematica} {\bf 29}, 217--234.

\end{document}